\documentclass[12pt]{article}
\usepackage{amsmath}
\usepackage{amsthm}
\usepackage{amsfonts}
\usepackage{amssymb}
\usepackage{latexsym}
\usepackage{euscript}
\usepackage{a4wide}

\def\Ta{\mathrm{Ta}}
\def\tors{\mathrm{tors}}

\renewcommand\mod{\ \mathrm{mod} \ }
\def\frak{\mathfrak}

\def\SL{\mathrm{SL}}

\def\GL{\mathrm{GL}}

\def\el{\ell}
\def\Qbar{\overline{\Q}}

\def\Gal{\mathrm{Gal}}
\def\F{\mathbf F}
\def\Fbar{\overline{\F}}

\def\rhobar{\overline{\rho}}

\def\qed{\hfill \square \ }
\def\Ot{\mathcal O}

\def\C{\mathbf C}
\def\I{I}

\def\Q{\mathbf Q}
\def\Z{\mathbf Z}
\def\R{\mathbf R}

\def\T{\mathbf T}

\def\el{\ell}
\def\GL{\mathrm{GL}}
\def\Pic{\mathrm{Pic}}
\def\Aut{\mathrm{Aut}}
\def\Sym{\mathrm{Sym}}

\def\et{\text{\'et}}

\def\ur{\mathrm{ur}}
\def\co{0}

\def\Hom{\mathop{\mathrm{Hom}}\nolimits}

\def\Tan{\mathop{\mathrm{Tan}}\nolimits}
\def\Spec{\mathop{\mathrm{Spec}}\nolimits}

\def\sdp{\medspace \times \kern -2.1pt 
 \vrule width0.4pt height 6pt depth -0.1pt \medspace}
\def\m{\frak{m}}
\def\ilim{\displaystyle \lim_{\longrightarrow}}

\def\f1{f'}

\def\Np{N\kern-.1em{p}}
\def\Npa{N'\kern-.1em{p}}
\def\Nm{N\kern-.1em{m}}
\def\Nma{N'\kern-.1em{m}}
\def\XYZ{XY\kern-.17em{Z}}

\def\iso{\cong}
\def\old{\mathrm{old}}

\swapnumbers

\newtheorem{theorem}[subsection]{Theorem}
\newtheorem{lemma}[subsection]{Lemma}

\newtheorem{cor}[subsection]{Corollary}

\newtheorem{prop}[subsection]{Proposition}
\newtheorem{remark}[subsection]{Remark}
\newtheorem{example}[subsection]{Example}

\def\XYZ{XY\kern-.17em{Z}}
\def\Kf{K_{\kern-.1em{f}}}

\def\f1{f'}

\begin{document}

\author{Frank Calegari\footnote{Supported in part by the American Institute
of Mathematics} \and Matthew Emerton\footnote{Supported in part by NSF
grant DMS-0401545}}
\title{Elliptic Curves of Odd Modular Degree}
\maketitle

\section{Introduction}

Let $E$ be an elliptic curve over $\Q$ of conductor $N$.
Since $E$ is modular \cite{bcdt},
there exists a surjective map $\pi: X_0(N) \rightarrow  E$ defined
over $\Q$.
There is a unique such map of minimal degree (up to composing with
automorphisms of $E$), and its degree $m_E$  is
known as the \emph{modular degree} of $E$.
In \cite{will}
the parity of $m_E$ was
determined for a very particular explicit
class of elliptic curves, namely, the
Neumann--Setzer curves, which have prime conductor and
a rational $2$-torsion point.  (See also the remark following
theorem~\ref{thm:NS} below.)
We study the question 
more generally for arbitrary elliptic curves $E/\Q$, 
and prove the following
theorem:

\begin{theorem} \label{theorem:one} 
If $E/\Q$ is an elliptic curve of odd modular degree then:

\begin{enumerate}

\item
the conductor $N$ of $E$ is divisible by at most two odd primes,
\item
$E$ is of even analytic rank,
and
\item
either
\subitem (a) $E$ has a rational point of order $2$ \emph{(}or equivalently,
admits a rational $2$-isogeny\emph{)},
\subitem
(b) $E$ has prime conductor and supersingular reduction at $2$, and $\Q(E[2])$
is totally complex \emph{(}equivalently, $E(\R)$ is connected\emph{)}, or
\subitem (c) $E$ has complex multiplication, and $N = 
27,32,49,\text{ or }243$.
\end{enumerate}
\end{theorem}

\begin{example}
{\em
The following 
examples of elliptic curves with odd modular degree should serve
to illustrate conditions~(3a), (3b) and~(3c).
The curve $X_0(15)$ has modular degree one and a rational
two torsion point, and thus satisfies condition~(3a).
Another example is given by the curve
$$y^2 + x y = x^3 - x^2 - 58x - 105$$
($2537E$ in Cremona's tables) of conductor $43 \cdot 59$ with modular
degree $445$ and torsion subgroup $\Z/4\Z$.
The curves $X_0(11)$ and $X_0(19)$ both have modular degree one and
satisfy condition~(3b).  An example of larger conductor
is given by
$$y^2 + y = x^3 + x^2 - 4 x - 10$$
of conductor $24859$ and
modular degree $3979$.
Finally, there are exactly four curves of odd modular degree with complex
multiplication,
namely $X_0(27)$, $X_0(32)$, $X_0(49)$ (all of modular degree one) and
$$y^2 + y = x^3 + 2$$
of modular degree $9$, conductor $243$ and $j$-invariant $0$.}
\end{example}

\begin{remark}
\label{remark:optimal}
{\em
Each of the conditions appearing in theorem~\ref{theorem:main}
is invariant under isogeny, other than
the condition that $E(\R)$ be connected, which however
is invariant under isogenies of odd degree.
Since the modular parameterization of $E$ factors through the optimal
member of the isogeny class of $E$
(that is, the member of its isogeny class having minimal modular degree;
in older terminology, a strong Weil curve),
it is therefore no loss of generality in the proof of theorem~\ref{theorem:one}
to assume that $E$ is optimal.}
\end{remark}

\begin{remark}
{\em
Cremona and Watkins have computed the modular degree of all optimal elliptic
curves of conductor $\leq$ 25,000 
\cite{CW}
and these computations suggest that there may be even stronger
limitations on the conductor of a curve of odd modular degree
than those imposed by theorem~\ref{theorem:one}.  Indeed, in the
range of Watkins' computations,
every curve of
odd modular degree has conductor divisible by at most two primes,  
and the conductor always has one of the following forms:
$2 p$, $4 p$, or $p q$, where $p$ and $q$ are odd primes.}
\end{remark}

\begin{remark}
{\em
The statement of the theorem regarding the analytic rank
of $E$ is consistent
with the rank conjecture of Watkins \cite[conj.~4.1]{watkins} (and
the conjecture of Birch and Swinnerton-Dyer).}
\end{remark}

The following result, conjectured by Watkins
(\cite[conj.~4.3]{will}, \cite[conj.~4.2]{watkins}),
is a simple consequence of theorem~\ref{theorem:one} 
(see lemma~\ref{lemma:simple}):

\begin{theorem} \label{theorem:Wat} Let $E/\Q$ be an elliptic curve
of prime conductor $N$, and suppose that $E$ is neither
a Neumann--Setzer curve, nor $X_0(17)$
\emph{(}equivalently,
$E$ does not have a rational $2$-torsion point\emph{)}.
If $m_E$ is odd, then $N \equiv 3 \mod 8$.
\end{theorem}

\medskip

%

One technique for proving that an elliptic curve $E$ has
even modular degree is to show that the map $\pi$ factors through
$X_0(N)/w$ for some non-trivial Atkin--Lehner involution $w$.
We use this approach in section~2 to prove theorem~\ref{theorem:composite-main},
which in turn implies parts~(1) and~(2) of
theorem~\ref{theorem:one}, and shows that~(3a) holds if
$N$ is divisible by at least two primes.
It remains to prove~(3) in the 
case when $N$ is a prime power.
The most difficult case to handle is
when $N$ is actually prime,
and in this case we deduce theorem~\ref{theorem:one} 
from the following result, proved in section~3.

\begin{theorem} \label{theorem:main} Let $N$ be prime,
let $\T$ denote the Hecke algebra over $\Z$ acting on
weight two cuspforms on $\Gamma_0(N)$, and let
$\m$ be a maximal ideal of $\T$ such that $\T/\m = \F_2$,
and such that the associated semi-simple Galois representation
$\rhobar_{\m}: G_{\Q} \rightarrow \GL_2(\F_2)$
is irreducible.
If the completion $\T_{\frak{m}} = \Z_2$, then
\begin{enumerate}
\item $\rhobar$ is supersingular at $2$, and
\item $\rhobar$ is totally complex.
\end{enumerate}
\end{theorem}

The relevance of this result to theorem~\ref{theorem:one} is that,
since $N$ is prime in the context of theorem~\ref{theorem:main},
a result of Ribet~\cite{Ribet} shows (assuming, as we may, that $E$ is optimal)
that the modular degree of $E$
is even if and only if $2$ is a congruence prime for the newform
of level $N$ attached to $E$.

The proof of
theorem~\ref{theorem:main} is motivated by the following considerations:
If $p$ is an \emph{odd} prime and $\rhobar: G_{\Q}
\rightarrow \GL_2(\F_p)$ is a surjective modular representation,
then theorems of Wiles and Taylor--Wiles~\cite{Wiles, TW}
show that the universal minimal deformation ring $R_{\emptyset}$ attached
to $\rhobar$ is
isomorphic to the universal minimal modular deformation ring
$\T_{\emptyset}$ ($= \T_{\frak{m}}$, since $N$ is prime).
Since $\T_{\emptyset}$ is
a finite $W(\F_p) = \Z_p$ algebra with residue field $\F_p$,
it is exactly equal to $\Z_p$
if and only if it is an \emph{\'etale} $\Z_p$-algebra. On the other hand,
since $R_{\emptyset} \cong \T_{\emptyset}$, this is equivalent to
$R_{\emptyset}$ being \'etale over $\Z_p$, which is in turn equivalent
to the reduced Zariski cotangent space of $R_{\emptyset}$
being trivial.  Since by construction $R_{\emptyset}$ represents
the minimal deformation functor, its reduced Zariski cotangent 
space considered as a set has cardinality equal to the number of
minimal
deformations
$$\rho: G_{\Q} \rightarrow \GL_2(\F_p[x]/(x^2))$$
of $\rhobar$. Thus to prove that $\T_{\frak{m}} \ne \Z_p$ it suffices to show
that there exists a non-trivial minimal deformation of
$\rhobar$ to
$\GL_2(\F_p[x]/(x^2))$.

\medskip

In spirit, the proof of theorem~\ref{theorem:main} follows this strategy;
in other words, we determine whether or not $\T_{\m} = \Z_2$
by a calculation on tangent spaces.
A significant problem arises, however, since we are working
in the case $p=2$, whilst the method of Wiles and
Taylor--Wiles  applies only to $p > 2$. This is not a mere technical obstruction;
many phenomena can occur when $p=2$ that do not occur for odd $p$.
To name two such: the possible failure of $\T_{\m}$ to be
Gorenstein and the consequent failure of
multiplicity one \cite{Kilford}, and the fact that
$\rhobar$ can arise from a totally real extension of $\Q$.
Calculations
in the second case suggest that the Taylor--Wiles strategy for proving
$R = \T$ in the minimal case
will not work without some significant new idea, since the numerical
coincidences that occur for odd $p$ whilst balancing
the Selmer and dual Selmer groups in the Greenberg--Wiles product formula
(see for example the remarks of de Shalit \cite{Shalit}, top of p.~442) 
do not occur in the case $p=2$.
Mark Dickinson~\cite{Dick} has proved an $R=\T$
theorem for $p=2$; however, his result requires many non-trivial hypotheses,
and indeed does not apply to \emph{any} of the representations considered
in theorem~\ref{theorem:main}.  (Its main application to date has been to
representations with image $\SL_2(\F_4) \simeq A_5$.) The issue here is
that the Taylor--Wiles
auxiliary prime arguments fail when $p=2$ and the image
of $\rhobar$ is dihedral.

Thus, instead of appealing to any general modularity results,
we show that $\T_{\m}$ is bigger than $\Z_2$ by explicitly constructing
(in certain situations) 
non-trivial deformations of $\rhobar_{\m}$ to $\F_2[x]/(x^2)$
that are demonstrably modular (and hence contribute to the 
reduced cotangent space of $\T_{\m}$).  The most difficult point
is to show that these deformations are modular of the correct (minimal)
level.  We prove this via a level-lowering result for modular forms
with values in Artinian $\Z_2$-algebras (theorem~\ref{thm:level-lowering}
below).   This level lowering result may be of independent interest;
for example, it provides evidence that an $R = \T$ theorem should hold
for those $\rhobar$ of characteristic two to which it applies.

The proof of~(3) when $N$ is a prime power (but is not actually prime)
is given in section~4.  In section~5 we make some concluding
remarks.

Let us close this introduction by pointing out that recently
Dummigan \cite{Dumm} has provided a heuristic explanation for
Watkins' rank conjecture that also relies on a hypothetical
$R = \T$ theorem for the residual Galois representation $\rhobar$
arising from the $2$-torsion on an elliptic curve $E$:
he uses the symmetric square map from $\rhobar$ to $\Sym^2 \rhobar$
to lift elements from $2$-Selmer group of $E$
to the tangent space to the deformation ring of $\rhobar$.
He also shows that the resulting tangent space elements 
can {\em never} be ``trapped'' (in the words of \cite[p.~450]{Wiles})
by the Taylor--Wiles method of introducing auxiliary primes.

Altogether, the experimental work of Watkins on the parity
of modular degrees, taken together with
the results of this paper and of \cite{Dumm},
suggests the validity of an $R = \T$ theorem for (at least certain)
residual Galois representations arising from the $2$-torsion on elliptic
curves, a theorem whose proof, however,
seems out of reach of the current techniques.

\section{$N$ composite with at least two distinct prime factors}
\label{section:composite}

In this section we prove the following theorem.

\begin{theorem} \label{theorem:composite-main}
If $E$ is an elliptic curve of odd modular degree then
the conductor $N$ of $E$ is divisible by at most two odd primes,
and $E$ is of even analytic rank.
Furthermore, if $N$ is divisible by at least two primes,
then $E$ contains a rational $2$-torsion point.
\end{theorem}

We begin with a preliminary lemma.  Let $E$ be an elliptic
curve over a field $k$;
let $O$ denote the origin of $E$.
We let $A$ denote the group of 
automorphisms of $E$ as a curve over $k$ (i.e.~$k$-rational
automorphisms of $E$ that do not necessarily fix $O$),
and suppose that $W$ is a finite elementary abelian $2$-subgroup of $A$.

\begin{lemma}
\label{lemma:bound}
The order of $W$ divides twice the order of $E[2](k)$.
\end{lemma}

\begin{Proof}
Let $A_0$ denote the subgroup of $A$ consisting
of automorphisms of $E$ as an elliptic curve over $k$
(i.e.~$k$-rational automorphisms of $E$ that do fix $O$).
The action of $E(k)$ on $E$ via translation realizes
$E(k)$ as a normal subgroup of $A$ which
has trivial intersection with $A_0$,
and which together with $A_0$ generates $A$. 
Thus $A$ sits in
the split short exact sequence of groups
\begin{equation}
\label{equation:automorphism ses}
 0 \rightarrow E(k) \rightarrow A \rightarrow A_0 \rightarrow 1 .
\end{equation}
(This is of course well known. The surjection $A \rightarrow A_0$
may also be regarded as the map $A = \Aut(E) \rightarrow \Aut(\Pic^0(E))$
induced by the functoriality of the formation of Picard varieties
--- the target being the group of automorphisms of
$\Pic^0(E)$ as a group variety
--- once we identify $E$ and $\Pic^0(E)$ as group varieties in the usual way.)

The short exact sequence~(\ref{equation:automorphism ses})
induces a short exact sequence
$$ 0 \rightarrow W \cap E(k) \rightarrow W \rightarrow W_0 \rightarrow 1 , $$
where $W_0$ denotes the projection of $W$ onto $A_0$.
The known structure of $A_0$ 
shows that $W_0$ is either trivial or of order $2$.
Since $W\cap E(k) \subset E[2](k),$
the lemma follows.
$\qed$
\end{Proof}

\

\noindent
{\bf Proof of theorem~\ref{theorem:composite-main}}.
The discussion of remark~\ref{remark:optimal}
shows that it suffices to prove the
theorem under the additional assumption that $E$ is an optimal quotient
of $X_0(N)$.

Let $W$ denote the group of automorphisms of $X_0(N)$ generated
by the Atkin--Lehner involutions; $W$ is an elementary abelian $2$-group
of rank equal to the number of primes dividing $N$.
Since $E$ is an optimal quotient of $X_0(N)$,
the action of $W$ on $X_0(N)$ descends to an action on $E$.
If $w \in W$ were to act trivially on $E$,
then the quotient map $X_0(N) \rightarrow E$ would factor
through $X_0(N)/w$, contradicting our assumption that $E$
is of odd modular degree.
Thus lemma~\ref{lemma:bound} shows that $W$ has order
at most $8$, and hence that $N$ is divisible by at most $3$ primes.
Furthermore, if $N$ is divisible by more than one prime,
then it shows that $E[2](\Q)$ is non-trivial.

Suppose now that $N$ is odd, so that $X_0(N)$ and $E$ both
have good reduction at $2$.  We may then apply the
argument of the preceding paragraph over $\F_2$,
and so conclude
from lemma~\ref{lemma:bound} that $W$ has order at most $4$.
Hence $N$ is divisible by at most two primes.

If $E$ is of odd analytic rank, and if $f_E$ denotes
the normalized newform of level $N$ attached to $E$,
then $w_N f_E = f_E$,
and so the automorphism of $E$ induced by $w_N$ has
trivial image in $A_0$.
Thus $w_N$ acts on $E$ via translation by an element
$P \in E(\Q)$.  Since $w_N$ interchanges the cusps
$0$ and $\infty$ on $X_0(N)$, we see that $P = \pi(0)
- \pi(\infty)$ (where $\pi: X_0(N) \rightarrow E$
is a modular parameterization of $E$, chosen so
that $\pi(\infty) = O$).

The assumption that $E$ has odd analytic rank also implies
that $L(f_E,1) = 0$.  Since this $L$-value can be computed
(up to a non-zero factor)
by integrating $f_E$ from $0$ to $\infty$ in the upper half-plane,
we conclude that $P = O,$ and thus that $w_N$ acts trivially
on $E$. Hence $\pi$ factors through the quotient
$X_0(N)/w_N$ of $X_0(N)$, and so must be of even modular degree,
a contradiction.
$\qed$

\section{$N$ prime}\label{section:prime}

\subsection{Reductions}

\begin{lemma} Theorem~\ref{theorem:main} implies part~(3)
of theorem~\ref{theorem:one} for $N$ prime.
\label{lemma:ribet}
\end{lemma}

\begin{Proof}
Suppose that $E$ is an elliptic curve of conductor $N$,
assumed to be optimal in its isogeny class.
Let $f_E$ be the associated Hecke eigenform of level $\Gamma_0(N)$
and weight $2$.
From a theorem of Ribet~\cite{Ribet},  $2 | m_E$  if and only if $f_E$ satisfies
a congruence with another cuspidal eigenform of level $N$. The set of
cuspidal eigenforms (in characteristic zero) congruent to $f$ is 
indexed by 
$\Hom(\T_{\m} \otimes \Q_2,\Qbar_2)$. Thus $f_E$ satisfies no
non-trivial congruences if and only if $\T_{\m} \otimes \Q_2 = \Q_2$,
or equivalently if and only if 
  $\T_{\m} = \Z_2$. $\qed$
\end{Proof}

\medskip 

The following lemma shows that it is possible to detect congruences
from modular forms modulo two.

\begin{lemma} \label{lemma:local}
Let $\frak{m}$ be the maximal ideal in $\T$
associated to $\rhobar$, and suppose there
exist two distinct non-zero $q$-expansions $f$, $g$ with
coefficients in $\F_2$ such that
$\frak{m}^2 f = \frak{m}^2 g = 0$.
Then  $\T_{\frak{m}} \ne \Z_2$.
\end{lemma}

\begin{Proof} Since $\rhobar$ is irreducible, $f$
and $g$ are both cuspidal. If $\T_{\frak{m}} = \Z_2$,
then $(\frak{m}^2,2) = \frak{m}$. Thus $f$ and $g$
are both $q$-expansions which are killed by $\frak{m}$.
By multiplicity one for $q$-expansions \cite[prop.~9.3]{eisenstein}
it follows that $f = g$. $\qed$
\end{Proof}

\medskip

Of use to us will be the following theorem of Grothendieck on Abelian
varieties with semistable reduction
\cite[Expos\'e IX, prop.~3.5]{Groth}:

\begin{theorem}[Grothendieck] \label{theorem:Groth} Let $A$ be an
Abelian variety over $\Q$ with semistable reduction at $\el$. Let
$\I_{\el} \subset \Gal(\Qbar/\Q)$ denote a choice of inertia group at
$\el$.
Then the action of $\I_{\ell}$ on the $p^n$-division
points of $A$ for $p \ne \el$
is rank two unipotent; i.e., as an endomorphism, for
$\sigma \in \I_{\el}$,
$$(\sigma - 1)^2 A[p^n] = 0.$$
In particular, $\I_{\el}$ acts through its maximal pro-$p$ quotient,
which is procyclic.
\end{theorem}

Shimura proved that a modular form $f$ of weight $2$ and
level $\Gamma_0(N)$ gives rise to a modular Abelian variety $A_f$ 
in such a way that the $p$-adic representations $\rho_f$ attached to $f$
arise from the torsion points of $A_f$. For prime $N$, these varieties are
semistable at $N$,
and so we may apply the theorem above to deduce that for $p = 2$,
such representations $\rho$ restricted to $\I_N$ factor
through a pro-cyclic $2$-group. For representations $\rhobar$ with image inside
$\GL_2(\F_2) \simeq S_3$, this means in particular that the order of inertia
at $N$ is either $1$ or~$2$.

\medskip

Let us now consider a Galois representation
$\rhobar: G_{\Q} \rightarrow \GL_2(\F_2)$, arising from
a cuspidal Hecke
eigenform of level $N$, whose image is not contained in a Borel subgroup.
(This is equivalent to $\rhobar$ being irreducible, and also to the image of
$\rhobar$ having order divisible by $3$.)
Let $L$ be the fixed field of the kernel
of $\rhobar$; the extension $L/\Q$ is unramified outside $2$ and $N$.
If $L/\Q$ is unramified at $N$ then $\rhobar$ has Serre conductor $1$,
contradicting a theorem of Tate \cite{Tate}. Thus by the discussion above
we see that inertia at $N$ factors through
a group of order $2$, that $L/\Q$ is an $S_3$-extension,
and (hence) that $\rhobar$ is absolutely irreducible.
Let $K/\Q$ be a cubic subfield
of $L$, and let $F$ be the quadratic extension inside $L$. 
Since $\rhobar$ is finite flat at $2$, it follows from Fontaine's
discriminant bounds \cite{Fontaine}
that the power of $2$ dividing the discriminant
of $F/\Q$ is at most $4$. Thus  $F/\Q$ 
must be $\Q(\sqrt{\pm N})$ (as it is ramified at $N$).  The extension $L/F$ is unramified at 
the prime above $N$ and is ramified at $2$ if and only if $\rhobar$
is supersingular.

\begin{lemma} \label{lemma:simple}
If $\rhobar$ is supersingular at two and not totally real
then $N \equiv 3 \mod 8$. In particular, theorem~\ref{theorem:one} implies
theorem~\ref{theorem:Wat}. \end{lemma}

\begin{Proof} By class field theory
the quadratic field $F/\Q$ admits a degree three extension
ramified precisely at $2$ only if $2$ is unramified and inert in $F$.
This occurs if and only if  $N \equiv 3 \mod 8$ and $F = \Q(\sqrt{-N})$, or
$N \equiv 5 \mod 8$ and $F = \Q(\sqrt{N})$. Moreover if $F = \Q(\sqrt{N})$ then
$K$ and $L$ are totally real. $\qed$ \end{Proof}

\medskip

We shall prove theorem~\ref{theorem:main}
by showing in the following subsections
that if $\rhobar$ satisfies at least one of the following conditions:
\begin{enumerate}
\item $\rhobar$ is totally real;
\item $\rhobar$ is unramified at $2$;
\item $\rhobar$ is ordinary, complex,  and ramified at $2$;
\end{enumerate}
then $\T_{\m} \neq \Z_2$.

\subsection{$\rhobar$ is totally real}
\label{section:real}

The theory of modular deformations is not well
understood when $\rhobar$ is totally real. Thus
our arguments in this section are geometric. We use the
following theorem, due to Merel~\cite[prop.~5]{Merel}.
(This interpretation of Merel's result is due to
Agashe~\cite[cor.~3.2.9]{Agashe}).

\begin{theorem} Let $N$ be prime. Then $J_0(N)(\R)$ is 
connected.
\label{theorem:merel}
\end{theorem}

If we let $g$ denote the dimension of $J:=J_0(N)$ it follows
that  $J(\R) \simeq (\R/\Z)^{g}$,  $J(\R)^{\tors} \simeq (\Q/\Z)^g$
and $J[2](\R) = (\Z/2\Z)^g$.

\medskip

Let $J[2^{\infty}] = \ilim J[2^m]$.
Then $J[2^{\infty}]$ is a $2$-divisible group over $\Q$
admitting an action of $\T_2$.

\medskip

Since $\T_2$ is
finite and flat over the 
complete local ring $\Z_2$
there exists a decomposition
$$\T_{2} = \prod \T_{\m},$$
where the product is taken 
over the  maximal ideals $\m$  of $\T$
of residue characteristic two. 
If $g({\m})$ denotes the 
rank of $\T_{\m}$ over $\Z_2$,
then
$$\sum_{\m} g({\m}) = \mathrm{rank}(\T_2/\Z_2) = g.$$

If $J[\m^{\infty}]:=J[2^{\infty}] \otimes_{\T_2} \T_{\m}$, then
$J[2^{\infty}] \simeq \prod J[\m^{\infty}]$
(compare~\cite{eisenstein} \S 7, p.~91).
From lemma~7.7 of~\cite{eisenstein}
we see that $\Ta_{\m}J \otimes \Q_2$ is free
of rank two over $\T_{\m} \otimes \Q_2$
(where $\Ta_{\m}J := \Hom(\Q_2/\Z_2,J[\m^{\infty}](\Qbar))$ is
the  $\m$-adic Tate module of $J$),
and thus that
\begin{equation}\label{eqn:m-div rank}
J[\m^{\infty}](\C)\iso (\Q_2/\Z_2)^{2 g({\m})}.
\end{equation}

Let $J[2]_{\m} := J[2]\otimes_{\T_2}\T_{\m}$ be
the $2$-torsion subgroup scheme of $J[\m^{\infty}]$.

\begin{lemma} For all maximal ideals
$\m$ of residue characteristic two there is an equality
$$\dim_{\Z/2\Z}(J[2]_{\m}(\R)) = g({\m})$$
\label{lemma:real}
\end{lemma} 

\begin{Proof} The isomorphism~(\ref{eqn:m-div rank})
induces an isomorphism
$J[2]_{\m}(\C) \iso (\Z/2\Z)^{2 g(\m)}.$
Let $\sigma \in  \Gal(\C/\R)$ denote complex conjugation.
Then $(\sigma - 1)^2 J[2]_{\m}(\C) = 0$. Thus
$J[2]_{\m}(\R)$ (which is the kernel of
$\sigma -1$) has dimension at least $g(\m)$.
If $\dim_{\Z/2\Z}(J_{\m}[2](\R)) > g(\m)$ for
some $\m$, then since
$$J[2](\R) = \prod J[2]_{\m}(\R),$$
and since (as was noted above) 
$\dim_{\Z/2\Z}(J[2](\R)) = g$,
we would deduce the inequality:
$$g = \sum \dim_{\Z/2\Z}(J[2]_{\m}(\R)) >
\sum_{\m} g({\m}) = g,$$
which is absurd.
$\qed$
\end{Proof}

\medskip

Now let $\rhobar$ be a totally real (absolutely)
irreducible continuous modular representation
of $\Gal(\Qbar/\Q)$ into $\GL_2(\F_2)$
of level $\Gamma_0(N)$, and let $\m$ be the corresponding
maximal ideal of $\T$.
The main result of \cite{BLR} shows that the $\Gal(\Qbar/\Q)$-representation
$J[\m](\Qbar)$ is a direct sum of copies of $\rhobar$.
Thus, since $\rhobar$ is totally real, we find that
$$\dim_{\Z/2\Z} J[\m](\R) = \dim_{\Z/2\Z} J[\m](\C) \geq \dim_{\Z/2\Z}\rhobar
= 2.$$
Combining this inequality with 
the inclusion $J[\m](\R) \subseteq J[2]_{\m}(\R)$ 
and lemma~\ref{lemma:real} we find that $g(\m) \geq 2$,
and thus that $\T_{\m} \neq \Z_2$.

\subsection{$\rhobar$ is unramified at $2$}

Suppose that $L/\Q$ is unramified at $2$. This forces $\rhobar$
to be ordinary.
By the theory of companion
forms~\cite{Gross} one expects that $\rhobar$ 
arises from a mod $2$ form of level $\Gamma_1(N)$ and weight $1$.
Although the results of~\cite{Gross} do not apply in
this case, Wiese~\cite{Wiese} explicitly constructs such
forms when the image of $\rhobar$ is
dihedral, as it is in our situation.
(In fact the only difficult point of
Wiese's construction is the case when $\rhobar$ is totally real,
and this case of theorem~\ref{theorem:main}
is already covered by section~\ref{section:real}.)
Let $f$ be this companion (Katz) modular form of
weight one modulo $2$. Let $A$ be the Hasse invariant modulo $2$,
which is a modular form of level one with $q$ expansion given by $1$.
Then $A  f = f$ and $g = f^2$ are two distinct $q$-expansions
modulo $2$ of weight $2$ and level $N$. Moreover, one sees
that $(T_{\ell} - a_{\ell}) f = (T_{\ell} - a_{\ell}) g = 0$
for all odd $\ell$, and that $(T_2 - a_2) f = 0$ and
$(T_2 - a_2)^2 g = 0$.
Thus 
$\frak{m}^2 f = \frak{m}^2 g = 0$ and
therefore by lemma~\ref{lemma:local},  $\T_{\m} \ne \Z_2$
and we are done. 

\subsection{$\rhobar$ is ordinary, complex, and ramified at $2$}
\label{ss:ord}

Suppose that $\rhobar$ is ordinary, complex, and ramified at $2$.
It follows that $F/\Q$ is complex and ramified at $2$, and thus
that $F = \Q(\sqrt{-N})$ for some $N \equiv 1 \mod 4$. Moreover,
the extension $L/F$ is unramified everywhere. Since $N \equiv 1
\mod 4$ it follows that $H: = L(\sqrt{-1})$ is also unramified
everywhere over $F$. The field $H$ is Galois over $\Q$, and
clearly
\begin{equation}\label{eqn:gal iso}
\Gal(H/\Q) \simeq S_3 \times \Z/2 \Z.
\end{equation}
We may embed $S_3 \times \Z/2\Z$ into $\GL_2(\F_2[x]/(x^2))$
by fixing an identification of $S_3$ with $\GL_2(\F_2)$,
and mapping a generator of $\Z/2\Z$ to the matrix
$$\left( \begin{matrix} 1 + x & 0 \\ 0 & 1 + x \end{matrix} \right).$$
Composing the isomorphism~(\ref{eqn:gal iso}) with this embedding
yields a representation:
$$\rho : \Gal(H/\Q) \hookrightarrow \GL_2(\F_2[x]/(x^2)).$$
The representation $\rho$ has trivial determinant 
(equivalently, determinant equal to the mod $2$ cyclotomic
character). We also claim that $\rho$ is finite flat at two.
To check this, it suffices to prove this over $\Z^{\ur}_2$.
The representation $\rho |_{\Gal(\Qbar_2/\Q_2^{\ur})}$ factors
through a group of order $2$, and one explicitly sees that
it arises as the generic fibre of the group scheme
$(D \oplus D)/\Z^{\ur}_2$, where $D$ is the non-trivial extension
of $\Z/2\Z$ by $\mu_2$ considered in~\cite{eisenstein} (Prop 4.2,
p.~58).
Thus one expects $\rho$ to 
arise from an $\F_2[x]/(x^2)$-valued modular form of weight two and level $N$,
corresponding to a surjective map of rings 
$\T_{\m} \rightarrow \F_2[x]/(x^2)$.  This would follow if
we knew that $\T_{\m}$ coincided with
the minimal deformation ring associated to
$\rhobar$. Rather than proving this, however, we shall use weight one
forms to explicitly construct 
a weight two modular form giving rise to $\rho$.

\medskip

Let $\chi_{4N}$ be the  character of conductor $4N$ associated to $F$.
Consider two faithful representations
$$\psi_1: \Gal(L/\Q) \iso S_3 \hookrightarrow \GL_2(\C), \qquad
\psi_2: \Gal(H/\Q) \iso S_3 \times \Z/2\Z \hookrightarrow \GL_2(\C).$$
Since $F/\Q$ is complex, these dihedral representations are odd
and therefore give
rise to weight one modular forms $h_1$, $h_2$ in $S_1(\Gamma_1(4N),\chi_{4N})$.

\begin{lemma} The modular forms $h_1$, $h_2$ are ordinary at $2$, have coefficients
in $\Z$, and are congruent modulo $2$. Let
$$g = \frac{(h_2 - h_1)}{2} \in \Z[[q]].$$
 Then $g \mod 2$ is non-zero, and the form $h = h_1 + x g
\in S_1(\Gamma_1(4N),\F_2[x]/(x^2))$
 is an eigenform
for all the Hecke operators, including $U_2$. The
  associated $\GL_2(\F_2[x]/(x^2))$ Galois representation 
attached to $h$ is $\rho$.
\end{lemma}

\begin{Proof} The modular forms are both ordinary at $2$ because the
representations $\psi_1$ and $\psi_2$ have non-trivial subspaces on
which inertia at two is trivial (since $I_2$ acts through a group of order $2$).
They both have coefficients in $\Z$, since $2 \cos(\pi/3) \in \Z$. The
congruence $h_1 \equiv h_2$ follows from the fact that both are
ordinary-at-$2$
Hecke eigenforms, and that $a(h_1,\ell) = a(h_2,\ell) 
\chi_4(\ell)$ for all
odd primes $\ell$, where $\chi_4$ is the character of conductor $4$.
From this one also sees that $g$ is non-trivial modulo two.
The claims about $h$ follow formally from the fact that
 $h_1$ and
$h_2$ are Hecke eigenforms that are congruent modulo $2$,
and the definition of $\rho$.  $\qed$
\end{Proof}

\medskip

Now that we have constructed the weight one form $h$ of level $4N$
giving rise to $\rho$,
we would like to construct a corresponding
weight two form of level $N$. Multiplying $h$ by the Eisenstein
series in $M_1(\Gamma_1(4N),\chi_{4N})$, we see that $h$ is the $q$-expansion
of a modular form in $S_2(\Gamma_0(4N),\F_2[x]/(x^2))$. Since $h$ is ordinary
and is a $U_2$ eigenform, we may apply the $U_2$ operator to deduce that
$h \in S_2(\Gamma_0(2 N),\F_2[x]/(x^2))$. 
Applying theorem~\ref{thm:level-lowering} (proved in the following 
subsection) 
we then deduce that in fact $h \in S_2(\Gamma_0(N), \F_2[x]/(x^2))$,
and (thus) that
$g \in S_2(\Gamma_0(N),\F_2)$. An easy calculation shows
that if $f = h_1 \equiv h_2  \equiv \sum a_n q^n$ then
$ \frak{m} f = 0$, whilst $\frak{m} g \subset \F_2 f$, and 
so also $\frak{m}^2 g = 0$.   
Thus $\T_{\frak{m}} \ne \Z_2$,
by lemma~\ref{lemma:local}.

\subsection{Level-lowering for modular deformations}
\label{ss:level-lowering}

The goal of this section is to prove a level-lowering result
for modular forms with coefficients in Artinian rings that
strengthens the case $p=2$ of~\cite[thm.~2.8]{Edixhoven}
(which in turn extends a level lowering result proved by  
Mazur \cite[thm.~6.1]{level-lowering} in the odd prime case).

We first establish 
a version of the multiplicity one theorem \cite[thm.~2.1]{Wiles} for $p=2$.
Under the additional assumption that $\rhobar$ is not finite
at $2$, this theorem was proved in \cite[\S 2]{buzz} (as was the corresponding
result for odd level).  Thus the key point in our theorem is that
$\rhobar$ is allowed to be finite at $2$, even though the level is taken to
be even.

\begin{theorem}\label{thm:mult-one}
Let $N$ be an odd natural number,
and let $\T$ denote the full $\Z$-algebra of Hecke operators
acting on weight two cuspforms of level $2 N$.
If $\m$ is a maximal ideal in $\T$ whose residue field $k$
is of characteristic $2$, and for which the associated
residual Galois representation
$$\rhobar: \Gal(\Qbar/\Q) \rightarrow \GL_2(k)$$
is \emph{(}absolutely\emph{)} irreducible, ordinary, and ramified at $2$,
then
$\Ta_{\m}J_0(2 N)$
\emph{(}the $\m$-adic Tate module of $J_0(2 N)$\emph{)} is free of rank
two over the completion $\T_{\m}$.
\end{theorem}

To be clear, the condition ``ordinary at $2$'' means that
the image of a decomposition group at $2$ under $\rhobar$
lies in a Borel subgroup of $\GL_2(k)$.  Since $k$ is of
characteristic $2$, we see that (for an appropriate choice
of basis) the restriction of $\rhobar$ to an inertia group
at $2$ may be written in the form
$$\rhobar_{| I_2} = \left( \begin{matrix} 1 & * \\ 0 & 1 \end{matrix}
\right ).$$
The assumption that $\rhobar$ is ramified at $2$ then implies
that $*$ is not identically zero.  Thus the representation
space of $\rhobar$ has a unique line invariant under $I_2$, 
and so $\rhobar$ is irreducible if and only if it is absolutely
irreducible.

\begin{lemma}\label{lem:unique}
Let $k$ be a finite field of characteristic $2$.
If $\rhobar: \Gal(\Qbar_2/\Q_2) \rightarrow \GL_2(k)$
is a continuous representation 
that is finite, ordinary, and ramified at $2$,
then $\rhobar$ has a unique finite flat prolongation
over $\Z_2$ (up to unique isomorphism).  Furthermore,
this prolongation is an extension of a rank one \'etale 
$k$-vector space scheme by a rank one multiplicative $k$-vector
space scheme.
\end{lemma}

\begin{Proof}
%
Any finite flat group scheme that prolongs an unramified
continuous representation of $\Gal(\Qbar_2/\Q_2)$ 
on a one-dimensional $k$-vector space is either \'etale or multiplicative.
Thus there are {\it a priori} four possible structures
for a finite flat prolongation of $\rhobar$:
\'etale extended by \'etale; multiplicative extended by multiplicative;
multiplicative extended by \'etale; or \'etale extended by multiplicative.
However, all but the last possibility necessarily gives rise to an
unramified generic fibre (note that any extension of multiplicative by
\'etale must split, by a consideration of the connected \'etale sequence).
Thus, since $\rhobar$ is assumed ramified, we see that any prolongation
of $\rhobar$ must be an extension of a rank one \'etale $k$-vector space scheme
scheme by a rank one multiplicative $k$-vector space scheme.

To see that such a prolongation is unique,
consider the minimal and maximal prolongations $M$ and $M'$ of
$\rhobar$ to a finite flat group scheme 
\cite[cor.~2.2.3]{ray}.  The result of the preceding paragraph
shows that
the natural morphism $M\rightarrow M'$ necessarily induces an isomorphism
on the connected components of the identity, and on the corresponding
groups of connected components.  By the $5$-lemma, this morphism
is thus an isomorphism, and the lemma follows.
$\qed$
\end{Proof}

\

We now show that certain results of Mazur \cite{mazur}
cited in the proof of \cite[thm.~2.1]{Wiles} extend to the case $p=2$.
We put ourselves in the context of \cite[\S 1]{mazur}, 
and use the notation introduced there.  Namely, 
we let $K$ denote a finite extension of $\Q_p$ for some prime $p$,
and let $\Ot$ denote the ring of integers of $K$.
If $A$ is an abelian variety over $K$, then we let $A_{/\Ot}$
denote the {\it connected component of the identity of} the N\'eron model of
$A$ over $\Spec \Ot$.
For any power $p^r$ of $p$, the $p^r$-torsion subgroup scheme
$A[p^r]_{/\Ot}$ of $A_{/\Ot}$ is then a quasi-finite flat group
scheme over $\Spec \Ot$; we let $FA[p^r]_{/\Ot}$ denote its
maximal finite flat subgroup scheme, 
and $A[p^r]^{\co}_{/\Ot}$ denote 
the maximal connected closed subgroup scheme
of $A[p^r]_{/\Ot}$. 
Since we took $A_{/\Ot}$ to be the connected component of the N\'eron
model of $A$, the inductive limit
$FA[p^{\infty}]_{/\Ot} := \ilim FA[p^r]_{/\Ot}$
is a $p$-divisible group,
and $A[p^{\infty}]^{\co}_{/\Ot} := \ilim A[p^r]^{\co}_{/\Ot}$
is the maximal connected $p$-divisible subgroup of $FA[p^{\infty}]_{/\Ot}$.

The following proposition is a variation on
\cite[prop.~1.3]{mazur}, in which we allow 
the ramification of $K$ over $\Q_p$ to be unrestricted,
at the expense of imposing more restrictive hypotheses on the
reduction of the abelian varieties
appearing in the exact sequence under consideration.

\begin{prop}\label{prop:exact}
Let $0 \rightarrow A \rightarrow B \rightarrow C \rightarrow 0$
be an exact sequence of abelian varieties over $K$ such that
$A$ has purely toric reduction, whilst $C$ has good reduction.
Then the induced sequence of $p$-divisible groups
$$ 0 \rightarrow A[p^{\infty}]^{\co}_{/\Ot} \rightarrow B[p^{\infty}]^{\co}_{/\Ot}
\rightarrow C[p^{\infty}]^{\co}_{/\Ot} \rightarrow 0$$
is a short exact sequence of $p$-divisible groups over $\Spec \Ot$.
Equivalently, for any power $p^r$ of $p$, the induced sequence
$$ 0 \rightarrow A[p^{r}]^{\co}_{/\Ot} \rightarrow B[p^{r}]^{\co}_{/\Ot}
\rightarrow C[p^{r}]^{\co}_{/\Ot} \rightarrow 0$$
is a short exact sequence of finite flat group schemes over $\Spec \Ot$.
\end{prop}

\begin{Proof}
Since $A$ has purely toric reduction, the group scheme
$A[p^r]^{\co}_{/\Ot}$ is of multiplicative type for each $r$.  Thus it
necessarily maps isomorphically onto its scheme theoretic
image in $B_{/\Ot}$, and thus the induced map
$A[p^{\infty}]^{\co}_{/\Ot} \rightarrow B[p^{\infty}]^{\co}_{/\Ot}$
is a closed embedding.

Let $C'\subset B$ be an abelian subvariety chosen so that the 
induced map $C' \rightarrow C$ is an isogeny.  Then $C'$ also
has good reduction, and so $C'[p^{\infty}]^{\co}_{/\Ot}
\rightarrow C[p^{\infty}]^{\co}_{/\Ot}$ is an epimorphism of $p$-divisible
groups over $\Spec \Ot$.  Thus the induced map
$B[p^{\infty}]^{\co}_{/\Ot} \rightarrow C[p^{\infty}]^{\co}_{/\Ot}$
is also an epimorphism of $p$-divisible groups. 
A consideration of generic fibres shows that the kernel of
this surjection coincides with the scheme-theoretic image 
of $A[p^{\infty}]^{\co}_{/\Ot}$ in $B[p^{\infty}]^{\co}_{/\Ot}$,
and so the proposition is proved.
$\qed$
\end{Proof}

\

\noindent
{\bf Proof of theorem~\ref{thm:mult-one}}.
We closely follow the method of proof of 
\cite[thm.~2.1~(ii)]{Wiles} in the case when ``$\Delta_{(p)}$ is trivial
$\mod \m$'' (in the terminology of that proof; see \cite[pp.~485--488]{Wiles}).
If we let $A$ denote the connected part of the kernel of the map
$J_0(2 N) \rightarrow J_0(N) \times J_0(N)$ induced by Albanese functoriality
applied to the two ``degeneracy maps'' from level $2 N$ to level $N$,
then $A$ is an abelian subvariety of $J_0(2 N)$ having purely toric reduction
at $2$, whilst the quotient $B$ of $J_0(2 N)$ by $A$ has good reduction at $2$.
From proposition~\ref{prop:exact} we obtain (for any $r\geq 1$)
the short exact sequence
$$ 0 \rightarrow A[2^r]^{\co}_{/\Z_2} \rightarrow J_0(2 N)[2^r]^{\co}_{/\Z_2}
\rightarrow B[2^r]^{\co}_{/\Z_2} \rightarrow 0$$
of connected finite flat group schemes over $\Spec \Z_2$.
By functoriality of the formation of this short exact sequence,
and since $A$ is a $\T$-invariant subvariety of $J_0(2 N)$,
we see that this is in fact a short exact sequence of $\T$-module schemes.
Localizing at $\m$ induces the corresponding short
exact sequence
\begin{equation}\label{eqn:ses-one}
0 \rightarrow A[2^r]^{\co}_{\m/\Z_2} \rightarrow J_0(2 N)[2^r]^{\co}_{\m/\Z_2}
\rightarrow B[2^r]^{\co}_{\m/\Z_2} \rightarrow 0.
\end{equation}
Passing to $\Qbar_2$-valued points induces a short exact sequence
of $\Gal(\Qbar_2/\Q_2)$-modules
\begin{equation}\label{eqn:ses-two}
0 \rightarrow A[2^r]^{\co}_{\m}(\Qbar_2) \rightarrow
J_0(2 N)[2^r]^{\co}_{\m}(\Qbar_2)
\rightarrow B[2^r]^{\co}_{\m}(\Qbar_2) \rightarrow 0,
\end{equation}
which is a subexact sequence of the short exact sequence
of $\Gal(\Qbar_2/\Q_2)$-modules
\begin{equation}\label{eqn:ses-three}
0 \rightarrow A[2^r]_{\m}(\Qbar_2) \rightarrow
J_0(2 N)[2^r]_{\m}(\Qbar_2)
\rightarrow B[2^r]_{\m}(\Qbar_2) \rightarrow 0.
\end{equation}

Let $A[2^r]_{\m}(\Qbar_2)^{\chi}$
(respectively
$J_0(2 N)[2^r]_{\m}(\Qbar_2)^{\chi}$,
respectively
$B[2^r]_{\m}(\Qbar_2)^{\chi}$)
denote the maximal
$\Gal(\Qbar_2/\Q_2)$-subrepresentation of 
$A[2^r]_{\m}(\Qbar_2)$
(respectively
$J_0(2 N)[2^r]_{\m}(\Qbar_2)$,
respectively
$B[2^r]_{\m}(\Qbar_2)$)
on which the inertia group acts through the $2$-adic cyclotomic
character $\chi$.
The short exact sequence~(\ref{eqn:ses-three}) induces
an exact sequence
\begin{equation}\label{eqn:ses-four}
0 \rightarrow A[2^r]_{\m}(\Qbar_2)^{\chi} \rightarrow
J_0(2 N)[2^r]_{\m}(\Qbar_2)^{\chi}
\rightarrow B[2^r]_{\m}(\Qbar_2)^{\chi}.
\end{equation}

\begin{lemma}\label{lem:mult}
Each of the groups schemes appearing in the exact
sequence~(\ref{eqn:ses-one}) is
of multiplicative type, and the exact sequences
(\ref{eqn:ses-two}) and~(\ref{eqn:ses-four}) coincide
(as subsequences of~(\ref{eqn:ses-three})). 
\end{lemma}

\begin{Proof}
We first remark that~(\ref{eqn:ses-three}) is the
exact sequence of $\T_{\m}[\Gal(\Qbar_2/\Q_2)]$-modules underlying the
the corresponding exact sequence of $\T_{\m}[\Gal(\Qbar/\Q)]$-modules
$$0 \rightarrow A[2^r]_{\m}(\Qbar) \rightarrow
J_0(2 N)[2^r]_{\m}(\Qbar)
\rightarrow B[2^r]_{\m}(\Qbar) \rightarrow 0.$$
Since $\rhobar$ is assumed irreducible as a $k[\Gal(\Qbar/\Q)]$-representation,
each of the modules
appearing in this exact sequence is a successive extension of copies
of $\rhobar$.  The same is thus true of each of the modules
appearing in the exact sequence~(\ref{eqn:ses-three}).

Since $A$ has purely toric reduction,
it is clear that 
$A[2^r]^{\co}_{\m/\Z_2}$ is of multiplicative type,  
and so
\begin{equation}\label{eqn:inclusion}
A[2^r]^{\co}_{\m}(\Qbar_2) \subset
A[2^r]_{\m}(\Qbar_2)^{\chi}.
\end{equation}
Fix a filtration  $0 = W_0 \subset W_1 \subset \cdots \subset W_n
= A[2^r]_{\m}(\Qbar_2)$ 
of $A[2^r]_{\m}(\Qbar_2)$ for which the successive quotients
$W_{i+1}/W_i$ are isomorphic to $\rhobar$.
Since $A$ has purely toric reduction
the quotient $A[2^r]_{/\Q_2}/A[2^r]^{\co}_{/\Q_2}$ is Cartier dual
to $\hat{A}[2^r]^{\co}_{/\Q_2}$ (where $\hat{A}$ is the dual abelian
variety to $A$), and so
$A[2^r]_{\frak m}(\Qbar_2)/A[2^r]^{\co}_{\m}(\Qbar_2)$
is an unramified $\Gal(\Qbar_2/\Q_2)$-representation.  
Since $\rhobar$ is assumed ramified at $2$,
this implies that
$$W_{i+1} \bigcap A[2^r]_{\m}^{\co}(\Qbar_2) \subsetneq 
W_{i} \bigcap A[2^r]_{\m}^{\co}(\Qbar_2)$$
for each $i \geq 0$.
Furthermore,
$$W_{i} \not\subset A[2^r]_{\m}(\Qbar_2)^{\chi}$$ for each $i > 0$,
because $\chi \mod 2$ is trivial.
Since $W_{i+1}/W_i \iso \rhobar$ is two dimensional over $k$
for each $i\geq 0,$ we conclude by induction on $i$ that
$$W_{i} \bigcap A[2^r]_{\m}^{\co}(\Qbar_2) =
W_{i} \bigcap A[2^r]_{\m}(\Qbar_2)^{\chi}$$
for each $i \geq 0$.  Taking $i = n$ then shows that
the inclusion~(\ref{eqn:inclusion}) is in fact an equality.

Since $B$ has good reduction at $2$, we have equality
$FB[2^r]_{/\Z_2} = B[2^r]_{/\Z_2}$.
As noted above, any Jordan--H\"older filtration
of the localization $B[2^r]_{\m}(\Qbar)$ 
as a $\T[\Gal(\Qbar/\Q)]$-module has
all its associated graded pieces isomorphic to $\rhobar$.
Taking scheme-theoretic closures of such a filtration in $B[2^r]_{/\Z_2},$
we obtain a filtration of the localization $B[2^r]_{\m/\Z_2}$
by finite flat closed subgroup schemes, whose associated graded pieces
are prolongations of $\rhobar$.  Now lemma~\ref{lem:unique} shows
that the connected component of any such finite flat prolongation is
multiplicative. 
Thus $B[2^r]^{\co}_{\m/\Z_2}$ is indeed multiplicative,
whilst $B[2^r]_{\m}(\Qbar_2)/B[2^r]_{\m}^{\co}(\Qbar_2)$
is an unramified $\Gal(\Qbar_2/\Q_2)$-module.  Arguing as in the
preceding paragraph gives the required equality
$$ B[2^r]_{\m}^{\co}(\Qbar_2) = B[2^r]_{\m}(\Qbar_2)^{\chi}.$$

Since any extension of multiplicative type groups is again
of multiplicative type, we see that
$J_0(2 N)[2^r]^{\co}_{\m/\Z_2}$ is also of multiplicative type,
and that the exact sequence~(\ref{eqn:ses-two})
is a subsequence of the exact sequence~(\ref{eqn:ses-four}).
We have furthermore shown that first and third non-trivial terms
of these two sequences coincide.  This implies that these
exact sequences do indeed coincide.
$\qed$
\end{Proof}

\

Specializing lemma~\ref{lem:mult} to the case $r=1$ 
shows that $J_0(2 N)[2]^{\co}_{\m}(\Qbar_2)$
is the maximal unramified
$\Gal(\Qbar_2/\Q_2)$-subrepresentation
of $J_0(2 N)[2]_{\m}(\Qbar_2)$
(since $\chi \mod 2$ is trivial).
Recall that there is a natural isomorphism
$\Tan(J_0(2 N)[2]^{\co}_{/\Fbar_2}) \iso \Tan(J_0(2 N)_{/\Fbar_2})$
(indeed, this is true with $J_0(2 N)_{/\Fbar_2}$ 
replaced by any group scheme over $\Fbar_2$),
and also a natural isomorphism
$\Tan(J_0(2 N)[2]^{\co}_{/\Fbar_2}) \iso
J_0(2 N)[2]^{\co}(\Qbar_2)\otimes_{\F_2} \Fbar_2$
(as follows from the discussion on \cite[p.~488]{Wiles}).
Localizing at $\m$, and taking into account \cite[lem.~2.2]{Wiles},
which is valid for $p=2$, we find that
$J_0(2 N)[2]^{\co}_{\m}(\Qbar_2)$ is a cyclic $\T_{\m}$-module,
and thus that the maximal unramified $\Gal(\Qbar_2/\Q_2)$-subrepresentation
of $J_0(2 N)[2]_{\m}(\Qbar_2)$ is a cyclic $\T_{\m}$-module.

Let $\rho_{\m}: \Gal(\Qbar/\Q) \rightarrow \GL_2(\T_{\m})$ denote
the Galois representation associated to $\m$ by \cite[thm.~3]{carayol}.
Carayol has proved \cite[thm.~4]{carayol} that there is an isomorphism
$\Ta_{\m}J_0(2 N) \iso J \otimes_{\T_{\m}} \rho_{\m}$
for some ideal $J$ in $\T_{\m}$,
and thus an isomorphism
$J_0(2 N)[2]_{\m}(\Qbar_2) \iso (J/2 J)\otimes_{\T_{\m}} \rho_{\m}$.
We conclude that $J/2 J$ is a cyclic $\T_{\m}$-module,
and hence that $J$ is a principal ideal in $\T_{\m}$.
The discussion of \cite[3.3.2]{carayol} shows that in fact $J \iso \T_{\m}$
and that $\Ta_{\m}J_0(2 N)$ is free of rank two over $\T_{\m}$, as claimed.
$\qed$

\

\begin{cor}\label{cor:gorenstein}
In the situation of theorem~\ref{thm:mult-one},
the completion $\T_{\m}$ is a Gorenstein $\Z_2$-algebra.
\end{cor}

\begin{Proof}
This follows from the theorem together
with the self-duality of the $\m$-adic Tate module
under the Weil pairing.
$\qed$
\end{Proof}

\

We now prove our level lowering result.
Let $A$ be an Artinian ring with
finite residue field $k$ of characteristic $2$,
and suppose given a continuous representation
$\rho: \Gal(\Qbar/\Q) \rightarrow \GL_2(A)$ that
is modular of level $\Gamma_0(2 N)$ for some odd natural number $N$,
in the sense that it arises from a Hecke eigenform
$h \in S_2(\Gamma_0(2 N),A)$. 
Let $\rhobar$ denote the residual representation attached to $\rho$
(so $\rhobar$ arises from the Hecke eigenform $\overline{h}
\in S_2(\Gamma_0(2 N), k)$ obtained by reducing $h$ modulo
the maximal ideal of $A$).

\begin{theorem}\label{thm:level-lowering}
If $\rho: \Gal(\Qbar/\Q) \rightarrow \GL_2(A)$ 
is a modular Galois representation
of level $\Gamma_0(2 N)$ as above, such that 
\begin{enumerate}
\item $\rhobar$ is \emph{(}absolutely\emph{)} irreducible,
\item $\rhobar$ is ordinary and ramified at $2$, and
\item $\rho$ is finite flat at $2$,
\end{enumerate}
then $\rho$ arises from an $A$-valued Hecke eigenform of level $N$.
\end{theorem}

\begin{Proof}
The Hecke eigenform $h$ corresponds to a ring homomorphism
$\phi: \T \rightarrow A$.    Since $A$ is local of residue
characteristic $2$, the map $\phi$ factors through the completion
$\T_{\m}$ of $\T$ at some maximal ideal $\m$
of residue characteristic $2$,  
and the residual representation $\rhobar$ is {\it the}
residual Galois representation attached to the maximal ideal $\m$.
We let $\rho_{\m}$ denote the Galois representation 
$$\rho_{\m}: \Gal(\Qbar/\Q) \rightarrow \GL_2(\T_{\m})$$
attached to $\m$ 
by \cite[thm.~3]{carayol}.
The Galois representation
$\rho$ attached to $h$ coincides with the pushforward of $\rho_{\m}$
via $\phi$.

Replacing $A$ by the image of $\phi$,
we may and do assume from now on that $\phi$ is surjective.
We let $I\subset \T_{\m}$ denote the kernel of $\phi$.
Since $A$ is Artinian, we may choose $r\geq 1$ so that $2^r \in I$.
Theorem~\ref{thm:mult-one} shows that the $\m$-adic Tate module
$\Ta_{\m} J_0(2 N)$ is isomorphic as a $\Gal(\Qbar/\Q)$-representation
to $\rho_{\m}$; thus $J_0(2 N)[2^r]_{\m}(\Qbar)$
is isomorphic to the reduction mod $2^r$ of $\rho_{\m}$.
Since $\T_{\m}$ is a Gorenstein $\Z_2$-algebra,
by corollary~\ref{cor:gorenstein},
we see that $\T_{\m}/2^r \T_{\m}$ is a Gorenstein
$\Z/2^n\Z$-algebra,
and thus that there is an isomorphism
$(\T_{\m}/2^r \T_{\m})[I]
\iso \Hom_{\Z/2^r}(\T_{\m}/I,\Z/2^r)$
of $\T_{\m}/I = A$-modules.
In particular, 
$J_0(2 N)[I](\Qbar/\Q) \subset J_0(2 N)[2^r]_{\m}(\Qbar)$
is a faithful $A$-module,
isomorphic as an $A[\Gal(\Qbar/\Q)]$-module to 
$\Hom_{\Z/2^r}(\T_{\m}/I,\Z/2^r)\otimes_A \rho.$
To simplify notation,
we will write
\begin{equation}\label{eqn:V}
V := J_0(2 N)[I](\Qbar) \iso
\Hom_{\Z/2^r}(\T_{\m}/I,\Z/2^r)\otimes_A \rho.
\end{equation}

By assumption, $\rho$ prolongs to a finite flat group scheme $\mathcal M$
over $\Spec \Z_2$. 
If we
fix a Jordan--H\"older filtration of $\rho$ as an $A[\Gal(\Qbar/\Q)]$-module,
then the associated graded pieces are each isomorphic to $\rhobar$,
and so lemma~\ref{lem:unique} and \cite[Prop.~2.5]{CE} together imply
that $\mathcal M$ is uniquely determined by $\rho$, whilst \cite[Lem.~2.4]{CE}
then implies that $\mathcal M$ is naturally an $A$-module scheme.
From~(\ref{eqn:V}) we see that $V$ also prolongs to a finite
flat $A$-module scheme
$$\mathcal V \iso \Hom_{\Z/2^r}(\T_{\m}/I,\Z/2^r)\otimes_A \mathcal M$$
over $\Z_2$.
Again,
lemma~\ref{lem:unique} and \cite[Prop.~2.5]{CE} show
that $\mathcal V$ is the unique finite flat prolongation of $V$.

Lemma~\ref{lem:unique} furthermore implies that $\mathcal M$ is
the extension of an \'etale $A$-module scheme $\mathcal M^{\et}$
by a multiplicative
$A$-module scheme $\mathcal M^{\co}$,
each of which is free of rank one as an $A$-module scheme.
Thus $\mathcal V$ is also an extension of an \'etale $A$-module scheme
$\mathcal V^{\et}$
by multiplicative $A$-module scheme $\mathcal V^{\co}$,
each of which is faithful as an $A$-module scheme.  
Let $V^{\et}$ and $V^{\co}$ denote the generic fibres of these schemes.

We write $\mathcal J$ to denote the
N\'eron model of $J_0(2 N)$ over $\Spec \Z_2$.
For a scheme over $\Z_2$, use the subscript ``$s$'' to
denote its special fibre over $\F_2$.
The special fibre $\mathcal J_s$ admits the following filtration
by $\T$-invariant closed subgroups:
$$0 \subset T \subset \mathcal J_s^{\co} \subset \mathcal J_s,$$
where $T$ is the maximal torus contained in $\mathcal J_s,$
and $\mathcal J_s^{\co}$ is the connected component of the identity
of $\mathcal J_s$.  The quotient $\mathcal J_s^{\co}/T$ is an abelian variety
on which $\T$ acts through its quotient $\T_{\old}$
(where $\T_{\old}$ denotes the quotient of $\T$ that acts faithfully
on the space of $2$-old forms of level $2 N$).
The connected component group
$\Phi:= \mathcal J_s/\mathcal J_s^{\co}$ is Eisenstein
\cite[Thm.~3.12]{level-lowering}.

The following lemma provides an analogue in our situation of
\cite[lem.~6.2]{level-lowering} (and generalizes one step
of the argument in the proof of \cite[thm.~2.8]{Edixhoven}). 

\begin{lemma}\label{lem:neron}
The Zariski closure of $V$ in 
$\mathcal J$
is a finite flat $A$-module scheme over $\Z_2$
(which is thus isomorphic to $\mathcal V$).
\end{lemma}

\begin{Proof}
Since $\mathcal V^{\co}$ is a multiplicative type group scheme,
inertia at $2$ acts on $V^{\co}(\Qbar_2)$ through
the cyclotomic character.  It follows from lemma~\ref{lem:mult}
that $V^{\co}$ is contained in the generic fibre of
$J_0(2 N)[2^r]_{\m/\Z_2}^{\co}$, and thus that the Zariski
closure of $V^{\co}$ in $\mathcal J$ is indeed finite flat,
and in fact of multiplicative type.  Thus it coincides
with $\mathcal V^{\co}$, and so we see that the embedding
of $V^{\co}$ in $J_0(2 N)$ prolongs to an embedding of
$\mathcal V^{\co}$ in $\mathcal J$.  Since the quotient
$\mathcal V^{\et} = \mathcal V/\mathcal V^{\co}$ is \'etale,
lemma~5.9.2 of \cite[Expos\'e IX]{Groth} serves to complete
the proof of the lemma.
$\qed$
\end{Proof}

\

Lemma~\ref{lem:neron} allows us to regard $\mathcal V$
as a closed $\T$-submodule scheme of $\mathcal J$, and thus
to regard $\mathcal V_s$ as a closed $\T$-submodule scheme of $\mathcal J_s$.
Since $\rhobar$ is irreducible and $\Phi$ is Eisenstein,
we see that $\mathcal V_s$ is in fact
contained in $\mathcal J_s^{\co}$.  On the other hand, since $T$
is a torus, we see that $\mathcal V_s \bigcap T \subset \mathcal V_s^{\co}$.
Thus $\mathcal V^{\et}_s$ appears as a subquotient of $\mathcal J_s^{\co}/T$,
and in particular the $\T$-action on $\mathcal V^{\et}$ factors through
the quotient $\T_{\old}$ of $\T$.
Since $\mathcal V^{\et}_s$ is a faithful $A$-module scheme, we see that
the map $\phi:\T \rightarrow A$ factors through $\T_{\old}$,
completing the proof of the theorem.
$\qed$
\end{Proof}

\

We remark that the obvious analogue of theorem~\ref{thm:level-lowering}
in the case of odd residue characteristic
is also true.  The proof is similar but easier, relying
on the uniqueness results on finite flat models due to Raynaud \cite{ray}. 
Of course, in those cases when the $R = \T$ theorem of Wiles,
Taylor--Wiles, and Diamond
\cite{Diamond, TW, Wiles} applies, it is an immediate consequence
of that theorem.
(Thus theorem~\ref{thm:level-lowering} may be regarded as
evidence for an $R = \T$ theorem for those $\rhobar$ of residue
characteristic $2$ that satisfy its hypotheses.)

\section{$N$ a proper prime power}

There are only finitely many elliptic curves of conductor $2^k$ 
for all $k$, and we may explicitly determine which have odd modular 
degree. Therefore we assume that $E$ has conductor $N$, where $N = p^k$
with $k \ge 2$ and $p \ge 3$. Let $\chi$ be the unique quadratic character of
conductor $p$. Let $E'$ be the elliptic curve $E$ twisted by $\chi$. The
curve $E'$ also has conductor $N$, and moreover the associated modular
forms $f_E$ and $f_{E'}$ are congruent modulo $2$, since twisting by
quadratic characters preserves $E[2]$. Since $N$ is odd, any non-trivial
congruence modulo $2$ between $f_E$ and other
forms in $S_2(\Gamma_0(N))$ forces the modular degree $m_E$ to be even
\cite{Ribet}. Thus
we are done unless $f_E = f_{E'} = f_E \otimes \chi$. In particular, the
representations associated to $f_E$
must be induced from a quadratic field, and thus $E$ has complex multiplication.
(Alternatively, the equality $f_E = f_{E'}$ implies that
$E$ is isogenous to its twist and deduce this way that $E$ has CM.)
If $E$ has CM and prime power conductor then $E$ is one
of finitely many well known elliptic curves, for which we can directly
determine the modular degree by consulting current databases (for $N = 163^2$, we use the elliptic curve
database of Stein--Watkins, described in~\cite{ECDB}).

\section{Further remarks}

Certainly not every $E$ satisfying the conditions
of theorem~\ref{theorem:one} will actually have odd modular degree,
and one could try to refine this result by deducing additional
necessary conditions that $E$ must satisfy in order to have odd modular degree.
In this section we say a little about the related question of
whether or not $2$ is a congruence prime for the associated modular
form $f_E$, when $E$ satisfies either of conditions~(3a) or~(3b)
of the theorem.
 

For curves $E$ with a rational two torsion point, the 
modular form $f_E$ automatically satisfies a mod two congruence
with an Eisenstein series, and so detecting whether $f$ satisfies a
congruence with a cuspform is a more subtle 
phenomenon than in the non-Eisenstein situation.
One approach might be to relate the Hecke
algebra to an appropriate universal deformation ring (if the latter exists).
If $N$ is prime this can be done
\cite{CE}, and this enables
one to determine when $\T_{\m} = \Z_2$ for such representations.
The specific determination of when $\T_{\m} = \Z_2$, however,
was already achieved (for $N$ prime and $\rhobar$ reducible) 
by Merel in \cite{Merel}:

\begin{theorem}\label{thm:NS}
Let $N \equiv 1 \mod 8$ be prime, and let $\T_{\m}$
be the localization of the Eisenstein prime at $2$. Then
$\T_{\m} \ne \Z_2$ if and only if $N = u^2 + 16 v^2$ and
$v \equiv (N-1)/8 \mod 2$.
\end{theorem}

If $E$ is a Neumann--Setzer curve, then $N = u^2 + 64$ for some $u \in \Z$.
The result of Merel above then clearly implies that the optimal
Neumann--Setzer curve $E$ has odd modular degree if and only
if $N \not\equiv 1 \mod 16$. (An alternative proof of this fact,
relying on the results of \cite{eisenstein}, is given in
\cite{will}, thm.~2.1.)  If $E$ has composite conductor, then
one might try to generalize the results of
\cite{Merel} or \cite{CE} to this setting.

\medskip

Suppose now that
$E$ has prime conductor, that  $\rhobar$ is irreducible and supersingular, and
that $\Q(E[2])$ is totally complex.  Let $K$ and $L$ be
the extensions of $\Q$ attached to $E$ as in the discussion of
section~\ref{section:prime}.
If one had an $R = \T$ result of the type discussed in the introduction,
then to obtain further necessary conditions for $E$ to have odd
modular degree, it would suffice to establish sufficient conditions
for the existence of an
appropriate non-trivial minimal deformation $\rho: G_{\Q} \rightarrow
\GL_2(\F_2[x]/(x^2))$ lifting $\rhobar$. 
For representations
$\rhobar$ that were complex and ramified at $2$, but ordinary,
we constructed such a $\rho$ directly in subsection~\ref{ss:ord}
by considering a
quadratic genus field extension of $L$. When $E$ is supersingular, 
such deformations $\rho$ (when they exist) may be more subtle and can
not necessarily be constructed so directly. 
One can however prove the following result.

Recall that if $\rhobar$ is supersingular and totally
complex, then $K/\Q$ is totally ramified at $2$ and $K$ has
exactly two complex embeddings.

\begin{prop}
Suppose $\rhobar$ is supersingular at $2$ and totally complex.
If either:
\begin{enumerate}
\item the class number of $K$ is even, or
\item the fundamental unit $\epsilon$ of $\Ot_K$ satisfies
$v_{\pi}(\epsilon - 1) \ge 2$, where $\pi$ is the unique prime
above $2$ in $\Ot_K$,
\end{enumerate}

then there exists a non-trivial minimal deformation
$$\rho: \Gal(\Qbar/\Q) \rightarrow \SL_2(\F_2[x]/(x^2)).$$
\end{prop}

Note that $\SL_2(\F_2[x]/(x^2)) \iso S_4 \times \Z/2 \Z$.
The idea behind the proof of this result is to study $A_4$
extensions $H$ containing $L$ that are minimally ramified at $2$
and unramified at~$N$. Such extensions can be obtained by considering
the Galois closure over $\Q$
of certain quadratic extensions of $K$. One obtains
a suitable such extension either by considering an unramified extension of $K$,
in case~(1), or the extension $K(\sqrt{\epsilon_K})$, in case~(2).

\noindent \it Email addresses\rm:\tt \  fcale@math.harvard.edu\\
\tt \ emerton@math.northwestern.edu
\end{document}